\documentclass[reqno,12pt]{amsart}
\usepackage{amsmath, latexsym, amsfonts, amssymb, amsthm, amscd}
\usepackage{mathrsfs}
\usepackage[utf8]{inputenc} 
\usepackage[T1]{fontenc}
\usepackage[cal=boondoxo]{mathalfa}

\usepackage{color}
\usepackage{hyperref}
\usepackage{breakurl}

\usepackage{graphicx}

\newcommand{\R}{\mathbb{R}}

\newcommand{\N}{\mathbb{N}}

\newcommand{\gep}{\varepsilon}

\newfont{\indic}{bbmss12}

\def\d{\, \mathrm{d}}

\title{A brief and personal history of stochastic partial differential equations}
\author{Lorenzo Zambotti}
\address{
Laboratoire de Probabilit\'es, Statistique et Mod\'elisation \\ Sorbonne Université, Université de Paris, CNRS \\
4 Place Jussieu, 75005 Paris, France}
\email{zambotti\@@lpsm.paris}

\begin{document}

\begin{abstract}
We trace the evolution of the theory of stochastic partial differential equations from the foundation to its development, until
the recent solution of long-standing problems on well-posedness of the KPZ equation and the stochastic quantization in dimension three.

\end{abstract}

\maketitle


\section{Introduction}

In September 2017 I attended a meeting in Trento in honor of 
Luciano Tubaro, who was retiring. Mimmo Iannelli gave a humorous and affectionate talk whose title was {\it Abstract stochastic equations: when we used to
study in Rome's traffic jams}. He talked about the '70s, when he and Luciano were the first students of Giuseppe Da Prato's, who 
around 1975 proposed them to work on a brand new topic: stochastic partial differential equations. Since I was myself a 
PhD student of Da Prato's in the late '90s, on that day in Trento I was being told the 
story of the beginning of our scientific family. 

Then, a month later, I was at the Fields Institute in Toronto for a conference in honor of Martin Hairer, who had been awarded in 2014 a Fields medal {\it "for
his outstanding contributions to the theory of stochastic partial differential equations, and in particular for the creation of a theory of regularity structures for such equations"} (the official citation of the International Mathematical Union).

Within a few weeks I was therefore confronted with a vivid representation of the beginning of SPDEs and with a celebration of their culminating
point so far.  I realised that, because of Hairer's Fields medal, the mathematical community was suddenly aware of the existence of SPDEs, 
although very little was commonly known about them. 

For example, during his laudatio which introduced Hairer's talk at the 2014 International Congress of Mathematics in Korea,
Ofer Zeitouni felt the need to say to the audience {\it "I guess that many of you had never heard about stochastic partial differential equations"}. The other
three Fields medals in 2014 were awarded for work in, respectively, dynamical systems, Riemann surfaces and number theory. Certainly there was no need to introduce these topics to the 
mathematicians attending the ICM. However, after forty years of work, with thousands of published papers and hundreds of
contributors, SPDEs were still unknown to a large portion of the mathematical world. 

I decided to dedicate my talk in Toronto not just to Hairer's achievements, but to the whole community that had formed and nurtured him. 
In the last two years I have given several times this talk in different occasions. This special issue of DCDS gives me the opportunity to write down the
few thoughts I have to share about this topic, in the hope that someone else may continue this work and enrich this tale with other points of view.
I will make no claim to exhaustivity: the topic is vast and I know only a fraction of the literature.
I wish to explain the origin and the development of SPDEs from my personal point of view, and I apologise in advance for the aspects of this story
that I will fail to explain properly or even mention. I encourage anyone wishing 
to see this tale completed or told differently and better to do so and continue the work I am starting.

\section{The beginnings}

In principle, a Stochastic Partial Differential Equation (SPDE) is a Partial Differential Equations (PDE) 
which is perturbed by some random external force. This definition is however too general: if we have a PDE
with some random coefficients, where the randomness appears as a parameter and the equation can be set and
solved with classical analytic arguments, then one speaks rather of a \emph{random PDE}; this is the case for
example of a (deterministic) PDE with a random initial condition. 

A SPDE is, more precisely, a PDE which contains some stochastic process (or field) and cannot be defined with standard
analytic techniques; typically such equations require some form of stochastic integration. In most of the cases, 
the equation is a classical PDE perturbed by adding a random external forcing. One of the first examples is the
following \emph{stochastic heat equation with additive noise}
\begin{equation}\label{she}
\frac{\partial u}{\partial t} = \Delta u +\xi
\end{equation}
where $u=(u(t,x))_{t\geq 0,x\in\R^d}$ is the unknown solution and $\xi=(\xi(t,x))$ is the random external
force. Then one can add non-linearities and, in some cases, multiply the external
force by a coefficient which depends on the unknown solution, for example
\begin{equation}\label{she2}
\frac{\partial u}{\partial t} = \Delta u +f(u)+\sigma(u)\,\xi
\end{equation}
where $f,\sigma:\R\to\R$ are smooth. The product $\sigma(u)\, \xi$ is not always well-defined, since in many cases
of interest $\xi$ is a \emph{generalised function} and $u$ is not expected to be smooth; in this case one writes the equation in an integral form
and uses Itô integration to give a sense to the stochastic term.

The idea of associating PDEs and randomness was already present in the physics literature in the '50s and '60s, see for example \cite{spiegel52,lyon60,chen64,gibson67}. In the mathematical literature, several authors 
extended Itô's theory of stochastic differential equations (SDE) to a Hilbert space setting, see for example Dalecki\u{\i} \cite{dalecki66} and Gross \cite{gross67}.
In a paper published in 1969 \cite{zakai69}, Zakai wrote that the unnormalised conditional density in a filtering problem satisfies a linear SPDE. 

However, to my knowledge, the first papers which studied explicitly a SPDE as a problem in its own appeared in the '70s.
In 1970 Cabaña \cite{cabana70} considered a linear wave equation 
\[
\frac{\partial^2 u}{\partial t^2} + 2b\, \frac{\partial u}{\partial t}= \Delta u +\xi
\]
with a \emph{space-time white noise} $\xi$ and one-dimensional space variable $x$. This is a very important particular choice for the random external force: it is given by a
random generalised function $\xi$ which is \emph{Gaussian} and has very strong independence properties, namely
the "values" at different points in space-time are independent. 

In 1972 three papers were published on the topic: two in France (Bensoussan-Temam \cite{bt72} and Pardoux \cite{pardoux72}) and one in Canada (Dawson \cite{dawson72}). The French school was strongly influenced by
the PDE methods of the time, championed by Jacques-Louis Lions and his collaborators. Bensoussan and Temam \cite{bt72} considered
an evolution equation driven by a monotone non-linear operator $A_t$
\[
\frac{{\rm d}y}{{\rm d}t} + A_t(y)=\xi
\]
and with an external forcing $\xi$ which we can call
now \emph{white in time and coloured in space}; this means that values of the noise on points with different time-coordinate are independent, but there is a non-trivial correlation in space. In \cite{pardoux72} 
Pardoux considered a similar problem with multiplicative noise 
\[
\frac{{\rm d}y}{{\rm d}t} + A_t(y)=B_t(y)\,\xi
\]
where $B_t$ is a non-linear operator and the stochastic term is treated with the Itô integration theory. In 1975
Pardoux defended his PhD thesis written under the supervision of Bensoussan and Temam, which is considered the first extended work on the topic.

Dawson's paper \cite{dawson72} has a more probabilistic flavour. It treats the stochastic heat equations \eqref{she} and \eqref{she2} with
one-dimensional space variable $x$ and space-time white noise $\xi$; it shows that the solution
$u$ to the linear equation \eqref{she} is almost-surely continuous in $(t,x)$ (this is false in higher dimension, as we are going to see below); moreover, it introduces
the non-linear equation \eqref{she2} with the coefficient $\sigma(u)=\sqrt{u}$, which will soon become famous
as the equation of the Super-Brownian motion (for $f=0$). 

In the following years more and more researchers got interested in SPDEs. In particular, the Italian and
Russian schools were founded, respectively, in 1976 with Da Prato's first paper \cite{dpit76}
on the topic (together with his students Iannelli and Tubaro) and between 1974 and 1977 with
Rozovski\u{\i}'s papers \cite{rozovski74,rozovski75} and Krylov-Rozovski\u{\i}'s \cite{kr77}.

\section{The physical models}

In the '80s some theoretical physicists published a few very influential papers based on 
applications of SPDEs to several important physical problems: Parisi-Wu's \cite{pw81} and Jona Lasinio-Mitter's
\cite{jlm85} on the \emph{stochastic quantization}, 
and the Kardar-Parisi-Zhang model for the \emph{dynamical scaling of a growing interface} \cite{kpz86}. All these papers would be, thirty years later, an
important motivation for the theory of regularity structures, see below.

\subsection{The Stochastic Quantization}
The 1981 paper \cite{pw81} by Parisi and Wu proposed a dynamical approach to the construction of probability measures which arise in Euclidean Quantum Field Theory.
The difficulty with such measures is that they are supposed to be supported by spaces of \emph{distributions}  (generalised functions) on $\R^d$, which makes the definition
of \emph{non-linear} densities problematic. For example one would like to consider a measure on the space of distributions ${\mathcal D}'([0,1]^d)$ of the form
\[
\mu({\rm d}\phi)=\frac1Z \exp\left(-\int_{[0,1]^d} V(\phi(x))\d x\right) {\mathcal N}(0,(1-\Delta)^{-1})({\rm d}\phi)
\]
where ${\mathcal N}(0,(1-\Delta)^{-1})$ is a Gaussian measure with covariance operator $(1-\Delta)^{-1}$,
with $\Delta$ the Laplace operator on $[0,1]^d$ with suitable boundary conditions, and $V:\R\to\R$ is some potential. If $d>1$ then
${\mathcal N}(0,(1-\Delta)^{-1})$-a.s. $\phi$ is a distribution and not a function, and the non-linearity $V(\phi)$ is therefore ill-defined. Parisi-Wu
introduce a stochastic partial differential equation 
\begin{equation}\label{eq:pw}
\frac{\partial u}{\partial t}= \Delta u -u-\frac12\,V'(u) + \xi, \qquad x\in[0,1]^d
\end{equation}
which has $\mu$ as invariant measure, namely if $u(0,\cdot)$ has law $\mu$, then so has $u(t,\cdot)$ for all $t\geq 0$. This is an infinite-dimensional
analog of the classical \emph{Langevin dynamics}. By the ergodic theorem, for a
generic initial condition $u(0,\cdot)$, the distribution of $u(t,\cdot)$ converges to $\mu$ as $t\to+\infty$. Therefore one can use the stochastic dynamical
system $(u(t,\cdot))_{t\geq 0}$ in order to obtain useful information on $\mu$. 

We note however that, for $d>1$, the solution to \eqref{eq:pw} is expected to be again a distribution on space-time, at least this is the case for
the linear equation with $V'\equiv 0$. Therefore a rigorous study of this equation is also problematic, since $V'(u)$ is again ill-defined.

The first rigorous paper on the Parisi-Wu programme was by Jona Lasinio-Mitter \cite{jlm85}, where the authors chose the non-linearity $V(\phi)=\phi^4$
and the space dimension $d=2$, in order to construct the continuum $\phi^4_2$ model of Euclidean Quantum Field Theory \cite{simon74,gj87}, and called this 
equation the \emph{stochastic quantization}. Jona Lasinio-Mitter studied a modified version of equation \eqref{eq:pw} and obtained
probabilistically weak solutions via a Girsanov transformation; strong solutions to \eqref{eq:pw} were obtained in a later paper by Da Prato-Debussche \cite{dpd03}, see 
below. The case of space dimension $d=3$ remained however open until the inception of regularity structures.

\subsection{The KPZ equation}
The Kardar-Parisi-Zhang (KPZ) equation \cite{kpz86} is the following SPDE
\begin{equation}\label{eq:kpz}
\frac{\partial h}{\partial t}= \nu\Delta h +\lambda|\nabla h|^2 + \xi, \qquad x\in\R^d
\end{equation}
and describes the fluctuations around a deterministic profile of a randomly growing interface, where $\nabla$ is
the gradient with respect to the space variable $x$. 

From an analytic point of view, even if $d=1$ the KPZ equation is very problematic: if we consider the case $\lambda=0$ then
we are back to the stochastic heat equation with additive white noise \eqref{she}, for which it is known that the solution $u$ is not better than
Hölder-continuous in $(t,x)$ and certainly not differentiable; we expect $h$ in \eqref{eq:kpz} to have
at best the same regularity as $u$. In particular the gradient in space $\nabla h$ is defined only as
a distribution and the term $(\nabla h)^2$ is ill-defined. We restrict ourselves for simplicity to the case $\nu=\lambda=1/2$.

In the original KPZ paper \cite{kpz86} it was noticed that one can \emph{linearize} \eqref{eq:kpz} by means of the \emph{Cole-Hopf transformation}:
if we define $\psi=(\psi(t,x))_{t\geq 0,x\in\R}$ as the unique solution to the equation
\begin{equation}\label{eq:colehopf}
\frac{\partial \psi}{\partial t}= \frac12\frac{\partial^2 \psi}{\partial x^2} +\psi \, \xi, \qquad x\in\R,
\end{equation}
which is called the \emph{stochastic heat equation with multiplicative noise}, then $h:=\log\psi$ (formally) solves \eqref{eq:kpz}. 

In the first mathematical paper on KPZ, Bertini-Cancrini \cite{bc95} studied in 1995 the stochastic heat equation \eqref{eq:colehopf} in the Itô sense for $d=1$.
Since Mueller \cite{mueller91} had proved that a.s. $\psi(t,x)>0$ for all $t>0$ and $x\in\R$, then the Cole-Hopf solution $h=\log\psi$ is indeed well-defined.
Bertini-Giacomin \cite{bg97} proved in 1997 that
the stationary Cole-Hopf solution is the scaling limit of a particle system, the weakly-asymmetric simple exclusion process (WASEP); this celebrated
result was the first example of the \emph{KPZ universality class}, see below.

Since \eqref{eq:colehopf} is to be interpreted in the Itô sense, one can apply the Itô formula to $h=\log\psi$ and the result is, at least formally,
that $h$ solves
\begin{equation}\label{eq:kpz2}
\frac{\partial h}{\partial t}= \frac12\,\frac{\partial^2 h}{\partial x^2} +\frac12\left[(\partial_x h)^2-\infty\right] + \xi, \qquad x\in\R,
\end{equation}
which is almost \eqref{eq:kpz}, apart from the appearance of the famous infinite constant which is supposed
to \emph{renormalize} the ill-defined term $(\partial_x h)^2$. Making sense of this renormalization and
constructing a well-posedness theory for such equations were however
open problems for over 15 years until Hairer's breakthrough \cite{hairer13}, see below.

We note that the KPZ equation, and in particular its \emph{universality class}, has been one of the most fertile topics in probability theory of the last decade, with connections to particle systems, random matrices,
integrable probability, random polymers and much else. See the surveys by Quastel \cite{quastel12} and Corwin \cite{corwin16} for more details.

\subsection{Superprocesses}
SPDEs have also been applied to \emph{biological systems}, in particular in the context of the so-called \emph{superprocesses} introduced 
by Watanabe and Dawson in the '70s. Superprocesses are limits of discrete population models of the following type: particles evolve in a $\R^d$ space following 
some Markovian dynamic, typically Brownian motion, independently of each other; at random exponential times each particle dies and is replaced by a 
random number of identical particles, which become new elements of the population and behave as all other particles. We refer to 
the Saint-Flour lecture notes by Dawson \cite{dawson93} and Perkins \cite{perkins02} for pedagogical introductions to this topic.

The total number of members of the population which are alive at time $t\geq 0$ follows a standard branching process and is independent of the motion
of the particles. Therefore there are three situations, depending on the value $m$ of the average number of descendants that a particle has when
it dies: if $m>1$ the population grows at an exponential rate, if $m<1$ it dies after a finite and integrable time, if $m=1$ it dies after a finite but non-integrable time. The three situations are called, respectively, \emph{supercritical, subcritical} and \emph{critical}. 

The critical case, with Brownian spatial motion, has a scaling limit which is a Markov process with values in the space of measures on the state space $\R^d$; this process is called the \emph{super-Brownian motion}. 
If $d=1$, then Konno-Shiga \cite{ks88} proved in 1988 that a.s. this random measure has a continuous density $X_t(x)$ with respect to the Lebesgue measure ${\rm d}x$ on $\R$, and $(X_t(x))_{t\geq 0,x\in\R}$ solves the SPDE
\begin{equation}\label{eq:sbm}
\frac{\partial X}{\partial t}= \frac12\,\frac{\partial^2 X}{\partial x^2} + \sqrt{X}\,\xi.
\end{equation}
The diffusion coefficient of this equation, already introduced by Dawson in \cite{dawson72}, does not satisfy the usual Lipschitz condition and, indeed, \emph{pathwise uniqueness} for \eqref{eq:sbm} is still an open problem, see the papers by Mytnik-Perkins \cite{mp11} and Mueller-Mytnik-Perkins 
\cite{mmp14}. More precisely, the situation is the following: we consider the SPDE
\begin{equation}\label{eq:sigmaholder}
\frac{\partial X}{\partial t}= \frac12\,\frac{\partial^2 X}{\partial x^2} + \sigma(X)\,\xi,
\end{equation}
with $\sigma:\R\to\R$ a Hölder function with exponent $\gamma\in\,]0,1[$, namely $|\sigma(x)-\sigma(y)|\leq C|x-y|^\gamma$ and one looks
in general for solutions with values in $\R$, rather than in $\R_+$; in particular, for equation \eqref{eq:sbm} one would have $\sigma(u)=\sqrt{|u|}$. Then:
\begin{itemize}
\item if $\gamma>3/4$ 
we have pathwise uniqueness, namely if we have two solutions $(X^1,\xi)$ and $(X^2,\xi)$ to \eqref{eq:sigmaholder} driven by the same noise $\xi$ with
$X^1(0,\cdot)=X^2(0,\cdot)$ a.s., then $X^1\equiv X^2$ almost surely
\item if $\gamma<3/4$ then pathwise uniqueness fails in general and there are counterexamples
\item if $\sigma(0)=0$ and one is interested only in the class of \emph{non-negative} solutions, then it is not known whether pathwise uni\-queness holds
or fails in this class for $\gamma<3/4$. This leaves in particular the hope that the equation for super-Brownian motion \eqref{eq:sbm} may satisfy pathwise uniqueness.  However for the related equation of super-Brownian motion \emph{with immigration} the pathwise non-uniqueness was proved by Chen in \cite{chen15}.
\end{itemize}

If the state space $\R^d$ has dimension greater or equal to 2, then a.s. the measure $X_t({\rm d}x)$ is singular with respect to the Lebesgue measure
(see \cite{dh79}),
but the equation \eqref{eq:sbm} is still well-defined as a \emph{martingale problem}, since the diffusion coefficient $\sigma(x)=\sqrt{x}$ has the special
property that $\sigma^2(x)=x$ is linear. Remarkably, this martingale problem is well-posed and one can prove uniqueness in law of these superprocesses
using a technique called \emph{duality} due to Watanabe \cite{watanabe68}, see the cited paper by Konno-Shiga \cite{ks88}; duality can also be applied
to prove uniqueness for other processes, see the works of Shiga \cite{shiga81,shiga87} and Mytnik \cite{mytnik96}.

Finally, we mention that superprocesses are related to Le Gall's Brownian snake, see \cite{legall99}, which also plays a crucial role in the context of planar random maps, see e.g. Miermont's lecture notes \cite{miermont}.

\section{The theory}

During the '80s and the '90s, several monographs were published with the aim of presenting a systematic theory of SPDEs. 

The first major monograph was Walsh's Saint-Flour lecture notes \cite{walsh86}, which were published in 1986. In this course Walsh proposed
a general approach to SPDEs which has been very influential; his point of view has 
a very probabilistic flavour, since it consists in regarding the solution $u=u(t,x)$ of a (parabolic or hyperbolic) SPDE as a \emph{multi-parameter
process}, or more generally a \emph{multi-parameter random field}. The stochastic integration with respect to space-time white noise is developed according to this point of view, considering $t\mapsto \xi(t,\cdot)$ as a so-called \emph{martingale measure}, thus generalizing the Itô theory. We have used Walsh's notations for the equations numbered from \eqref{she} to \eqref{eq:sbm} above, and for others below.

In 1992 the first book by Da Prato-Zabczyk \cite{dpz1} was published. This monograph, also known as the \emph{red book} among Da Prato's students, is still the 
reference text for the so-called \emph{semigroup approach} to SPDEs. Da Prato-Zabczyk's point of view is to treat a SPDE as an-infinite dimensional SDE, and the solution $u=u(t,\cdot)$ as a function-valued process with a single parameter, the time $t$. The notations are different from those of Walsh;
for example the stochastic heat equation with additive space-time white noise \eqref{she} is written as
\[
\d X=AX\d t+\d W
\]
where $X_t=u(t,\cdot)\in L^2(\R)= H$, $A:D(A)\subset H\to H$ is the realization of $\partial^2_x$ in $H$, $(W_t)_{t\geq 0}$ is a \emph{cylindrical Wiener process}. The solution to this equation is called the \emph{stochastic convolution} and is written explicitly as
\[
X_t=e^{tA}X_0+\int_0^t e^{(t-s)A}\d W_s, \qquad t\geq 0.
\]
The general SPDE \eqref{she2} with non-linear coefficients is written as
\[
\d X=(AX+F(X))\d t+\Sigma(X)\d W
\]
where $F:D(F)\subseteq H\to H$ is some non-linear function and $\Sigma$ is a map from $H$ to the linear operators in $H$. This approach has a more 
functional-analytical flavour, and is based mainly on the study of the properties of the semigroup $(e^{tA})_{t\geq 0}$ generated by $A$ in $H$,
and their interplay with the properties of the cylindrical Wiener process $W$. This non-linear equation is usually written in its \emph{mild formulation}
\[
X_t=e^{tA}X_0+\int_0^t e^{(t-s)A}\,F(X_s)\d s+\int_0^t e^{(t-s)A}\,\Sigma(X_s)\d W_s.
\]

During the '90s there was also an important activity on \emph{infinite-dimensional analysis}, namely on elliptic and parabolic PDEs where the space-variable belongs to a Hilbert space. The connection with SPDEs is given by the notion of \emph{infinitesimal generator} which is associated with a 
Markov process with continuous paths. As for finite-dimensional diffusions, the transition semigroup of the solution to a SPDE solves a parabolic
equation, known as \emph{Kolmogorov equation}. One can find a systematic theory of these operators in the third book by Da Prato-Zabczyk \cite{dpz3}.
Much work was dedicated to existence and uniqueness of invariant measures, see the next section; the second Da Prato-Zabczyk book was
entirely dedicated to this topic \cite{dpz2}.

It can be recalled that Itô introduced his notion of stochastic differential equations in order to give
a probabilistic representation of the solution to Kolmogorov equations. Viceversa, if the Kolmogorov equation is well-posed, then it is
possible to construct the law of the associated Markov process. This allows to construct \emph{weak} (in the probabilistic sense) solutions, especially
in the form of \emph{martingale solutions}, see the 1979 monograph by Stroock-Varadhan \cite{sv79} on the theory for finite dimensional diffusions. 

The construction of the transition semigroup of a Markov process in a locally compact space can be done also with another analytical tool, a
\emph{Dirichlet Form}, for which a theory was developped in particular by Fukushima, see the monographs \cite{fukushima80,fot11}. The state space of a SPDE is
however always a function space, and therefore infinite-dimensional. The extension of Fukushima's theory to non locally compact spaces was 
a project of Albeverio-H{\o}egh-Krohn \cite{ah77} since the '70s and was finally obtained
by Ma-Röckner \cite{mr92}. Although Dirichlet forms allow to construct only weak solutions, they are a powerful tool in very singular situations, where
pathwise methods are often ineffective.

Another approach to SPDEs is given by Krylov's $L^p$-theory, see for example \cite{krylov94}.

\section{Ergodicity of Navier-Stokes}

The Navier-Stokes equation for the flow of an incompressible fluid is one of the most prominent PDEs and it is 
therefore not surprising that its stochastic version was among the first SPDEs to be studied,
starting from the 1973 paper \cite{bt73} by Bensoussan-Temam. The equation has the form (in Walsh's notation)
\[
\frac{\partial u}{\partial t}+(\nabla u)\cdot u= \nu\Delta u -\nabla p + \xi, \qquad {\rm div}\ u=0,
\]
where $u(t,x)\in\R^d$ denotes the value of the velocity of the fluid at time $t\geq 0$
and position $x\in\R^d$, $p(t,x)$ is the pressure, $\nu>0$ and $\xi$ is an external noise whose structure will be made 
precise below.

The statistical approach to hydrodynamics is based on the assumption that the fluid has a stationary state (invariant measure) on 
the phase space; by the ergodic theorem, the time average of an observable computed over the dynamics converges
for large time to the average of the observable with respect to the invariant measure. This ergodicity property must
however be proved, and in the case of the Stochastic Navier Stokes equation in 2D this has been a very active 
area of research, at least between the 1995 paper by Flandoli-Maslowski \cite{fm95} and the 2006 paper by
Hairer-Mattingly \cite{hm06}.

\subsection{Ellipticity versus hypoellipticity}
For stochastic differential equations in general, the choice of the external noise plays a very important role. 
In most of the literature on SPDEs, the 
space-time noise $\xi$ is realised as the following series
\[
\xi(t,x) = \sum_{k=1}^\infty \lambda_k\, e_k(x)\, \dot{B}_k(t), \qquad t\geq 0, \ x\in {\mathcal O}\subseteq \R^d,
\]
where $(\lambda_k)_k$ is a sequence of real numbers, $(e_k)_k$ an orthonormal basis of $L^2({\mathcal O},{\rm d}x)$ and $(B_k)_k$ an independent family of standard Brownian motions. If $\lambda_k=1$ for all $k$ then we have
space-time white noise, which has the property that for all $\varphi\in L^2({\mathcal O},{\rm d}x)$ the random 
variable
\[
\int_{[0,T]\times{\mathcal O}} \varphi(t,x)\, \xi(t,x)\d t\d x:=
\sum_{k=1}^\infty \langle \varphi,e_k\rangle_{L^2({\mathcal O},{\rm d}x)}\, {B}_k(T)
\]
has normal law ${\mathcal N}\left(0,T\,\|\varphi\|^2_{L^2({\mathcal O},{\rm d}x)}\right)$. 

In analogy with the finite-dimensional case, if $\lambda^2_k\geq \gep>0$ for all $k$, then we are in the 
\emph{elliptic} case. In finite dimension, we are in a degenerate case as soon as $\lambda_k=0$ for some
$k$; in infinite dimension, however, we can have $\lambda_k>0$ for all $k$ but $\lambda_k\to 0$ as $k\to+\infty$.
This situation is neither degenerate nor elliptic.

The paper by Flandoli-Maslowski proved for the first time ergodicity for a stochastic Navier-Stokes equation
in 2D, under the assumption that $\lambda_k>0$ for all $k$ but $\lambda_k\to 0$ as $k\to+\infty$ with two (different) power-law controls
from above and from below. This article sparked an intense activity and a heated debate which revolved around
the following question: what is the most relevant choice of the noise structure, which allows to prove ergodicity?

If, as in Flandoli-Maslowski \cite{fm95},  the noise is sufficiently non-dege\-nerate, namely if $\lambda_k>0$ and $\lambda_k\to 0$
not too fast as $k\to+\infty$, then it is often possible to prove
ergodicity using an argument due to Doob and based on two ingredients: the \emph{Strong-Feller property} and \emph{irreducibility}; the former means that
the transition semigroup of the dynamics maps bounded Borel functions on the state space into continuous functions,
the latter that all non-empty open sets of the state space are visited with positive probability at any positive time. The
Strong-Feller property is proved with ideas coming from Malliavin calculus, in particular on an integration by parts on the path space which is now known
as the Bismut-Elworthy-Li formula, see the paper by Elworthy-Li \cite{el94} and the monograph by Cerrai \cite{cerrai01};
irreducibility is based on control theory for PDEs. These
techniques were explored and applied to a number of examples in the second Da Prato-Zabczyk book \cite{dpz2} of 1996.

However, it soon appeared clear that it was possible to consider a degenerate noise and still obtain uniqueness of the invariant
measure. Here by degenerate we mean that $\lambda_k=0$ for all $k>N$, where $N$ is a deterministic integer. The main idea behind this
line of research was that, if the noise acted on a sufficiently large but finite number of \emph{modes} (i.e. the functions $e_k$), then the 
noise is elliptic on the modes which determine the long-time behavior of the dynamics:
we can call this the \emph{essentially elliptic} case. These results, together with exponential convergence to equilibrium, were proved independently (for Gaussian or for discrete noise) by three groups of authors during the same years: Mattingly \cite{matt99} and E-Mattingly-Sinai \cite{ems01}, Kuksin-Shirikyan \cite{ks00,ks01}, Bricmont-Kupiainen-Lefevere \cite{bkl01,bkl02}.

However in these works the number $N$ of randomly forced modes is not universal but depends on the parameters $\nu$ and $\sum_k \lambda^2_k$ of 
the equation. This was dramatically improved in the paper by Hairer-Mattingly \cite{hm06} published in 2006 in Annals of Mathematics, which proved that it is
enough to inject randomness only in \emph{four} well-chosen modes, then the non-linearity propagated the
randomness to the whole system for any $\nu>0$: the so-called \emph{hypoelliptic} case, for which it is possible to derive uniqueness of the invariant measure for the 2D stochastic Navier-Stokes. 
One of the main novelties in this paper was the notion of the \emph{asymptotic Strong-Feller property}, which could be
proved in the hypoelliptic case, while the standard Strong-Feller property requires much stronger non-degeneracy properties of the noise.

Let us mention here that the Malliavin Calculus, see e.g. Nualart's monograph \cite{nualart}, has played an important role for Navier-Stokes like for many other SPDEs.

\section{My SPDEs}

The results on the ergodicity of the stochastic Navier-Stokes equation seemed at the time to make SPDEs with degenerate noise particularly 
prominent. Now that singular SPDEs with space-time white noise and regularity structures have become so famous, this may seem even strange. 
In fact, since the very first papers that I have mentioned, see Cabaña \cite{cabana70} and Dawson \cite{dawson72}, the research activity on SPDEs with 
genuinely infinite-dimensional noise has always been intensive and most of the problems I have mentioned above concern space-time white noise. 

The case of degenerate noise is certainly more difficult if one wants to prove ergodicity, as we have seen. However, if the noise is spatially finite-dimensional, then the solution to the SPDE are typically smooth in space, although still Brownian-like in time. In the case of space-time white noise,
on the contrary, the solution are rather Brownian-like \emph{in space} if the space dimension is $d=1$, and even less regular in time; if $d>1$, as 
we have already seen, solutions are rather distributions. 

Therefore, SPDEs driven by space-time white noise are particularly strange objects: even the solutions to the simplest equation, as the
stochastic heat equation with additive space-time white noise,
are far too irregular for any of the derivatives which appear in \eqref{she} to make any sense as a function.
The KPZ equation \eqref{eq:kpz2} has almost an explicit solution given by the Cole-Hopf transform $h=\log\psi$, with $\psi$ solution to the
stochastic heat equation with multiplicative space-time white noise \eqref{eq:colehopf}; however the KPZ equation itself makes no sense as it is
written in \eqref{eq:kpz2}!

It is in this topic that I made my first steps as a researcher. I did my PhD at Scuola Normale in Pisa under the supervision of Giuseppe Da Prato (also
known as Beppe) from 1997 to 2001. Like Da Prato himself and many of his students, I started as an analyst but felt increasingly attracted by probability theory,
in particular stochastic calculus and SDEs. In the shelves of Beppe's office I found the Revuz-Yor monograph, which became one of my favourite mathematics books. 
I started to dream of unifying two worlds: the classical Itô theory of stochastic calculus based on martingales, and SPDEs.

Chapter 5 in the book by Revuz-Yor on local time and reflecting Brownian motion was one of the topics which most intrigued me. At that time Da Prato
was studying equations of the form
\begin{equation}\label{eq:sdi}
\d X\in (AX-\partial U(X))\d t+{\rm d} W
\end{equation}
with $U:H\to\R$ a \emph{convex} lower semi-continuous but not necessarily differentiable function. In the deterministic setting, this is a classical 
problem and the set $\partial U(x)$ is the \emph{subdifferential} at a point $x\in H$, namely the set of all directions $h\in H$ such that the affine subspace $U(x)+\{z\in H:\langle z,h\rangle =0\}$ lies below the graph of $U$. 
For a simple example, think of the function $\R\ni x\mapsto|x|\in\R_+$, which is convex and has as subdifferential the set $\{1\}$ for all $x>0$, the set
$\{-1\}$ for all $x<0$ and the set $[-1,1]$ for $x=0$. Then equation \eqref{eq:sdi} is rather a \emph{stochastic differential inclusion}, and if $U$ is differentiable at $x$ then $\partial U(x)=\{\nabla U(x)\}$. There is an extensive literature on this problem in the finite-dimensional case, see e.g. C\'epa \cite{cepa}, much less so in infinite dimension where many problems remain open.

The case of $U$ being equal to $0$ on a closed convex set $K\subseteq H$ and to $+\infty$ on $H\setminus K$ seemed to be outside the scope of Da Prato's techniques. I convinced myself that this case had to be related with reflection on the boundary of $K$, but I was unable to make this precise. Then Samy 
Tindel pointed out to me a 1992 paper by Nualart and Pardoux \cite{nupa} on the following SPDE with reflection at 0
\begin{equation}\label{eq:nupa}
\frac{\partial u}{\partial t} = \frac12\frac{\partial^2 u}{\partial x^2} +\xi +\eta, \qquad t\geq 0, \ x\in[0,1],
\end{equation}
where $\eta$ is a Radon measure on $]0,+\infty[\,\times\,]0,1[$, $u$ is \emph{continuous} and non-negative, and the support of $\eta$ is included
in the zero set $\{(t,x): u(t,x)=0\}$ of $u$, or equivalently
\begin{equation}\label{eq:nupa2}
u\geq 0, \qquad \eta\geq 0, \qquad \int_{]0,+\infty[\,\times\,]0,1[} u\d\eta=0.
\end{equation}
This is a \emph{stochastic obstacle problem}, the obstacle being the constant function equal to $0$, which can be formulated in the abstract setting of the
stochastic differential inclusion \eqref{eq:sdi}. 
Continuity of $(t,x)\mapsto u(t,x)$ here is essential in order to make sense of the condition \eqref{eq:nupa2}; in this setting the Walsh approach is clearly necessary, since continuity of $t\mapsto u(t,\cdot)$ in $L^2(0,1)$ would not be sufficient.
In higher space dimension, $u$ is not expected to be continuous and indeed it remains an open problem to define in this case a notion of solution to 
\eqref{eq:nupa}-\eqref{eq:nupa2}. We note also that this equation arises as the scaling limit of interesting microscopic models of random interfaces:
see Funaki-Olla \cite{fo01} and Etheridge-Labb\'e \cite{el15}.

The Nualart-Pardoux paper was motivated by stochastic analysis but it was an entirely deterministic work, which pushed the PDE techniques to 
cover a situation of minimal regularity for the solution; a probabilistic interpretation of this result remained elusive. This is what I tried to give
with the results of my PhD thesis. First I identified in \cite{lz01} the unique invariant measure of \eqref{eq:nupa}-\eqref{eq:nupa2} as the 3-d Bessel bridge (also known
as the normalized Brownian excursion), an important process which plays a key role in the study of Brownian motion and its excursion theory, see \cite{reyo}. Then I proved in \cite{lz02} an infinite-dimensional integration by parts with respect to the law of the 3-d Bessel bridge, which gave a powerful probabilistic
tool to study the reflection measure $\eta$ (it provides its \emph{Revuz measure}). 
Then I set out to study the fine properties of the solution, in particular of the contact set $\{(t,x): u(t,x)=0\}$ between the solution $u$ and the obstacle $0$, see \cite{lz04} and the paper \cite{dmz06} in collaboration with Dalang and Mueller.

In these papers I tried to realize my dream, by showing that solutions to SPDEs display very rich and new phenomena with respect to finite-dimensional
SDEs, and that it was possible to go much beyond results on existence and uniqueness. I found some interesting link between
classical stochastic processes arising in the study of Brownian motion and SPDEs. For a more recent account, see my Saint-Flour lecture notes
\cite{lz15}. 

However it does not seem that this point of view has been followed by many others. As we are going to see, the SPDE community would soon be 
heading in a very different direction.

\section{Rough paths and regularity structures}

In 1998 T. Lyons published a paper \cite{lyons98} on a new approach to stochastic integration. Lyons was an 
accomplished probabilist and an expert of stochastic analysis. Therefore it may seem puzzling that the aim of his most famous contribution to mathematics, 
the invention of \emph{rough paths}, is to give a deterministic theory of stochastic differential equations!

The classical Itô theory of stochastic calculus, see again \cite{reyo}, is a wonderful tool to study stochastic processes (more precisely continuous 
semimartingales). Not only does it allow to prove existence and uniqueness of solutions to stochastic differential equations, but it also allows to compute the
law of a great variety of random variables and stochastic processes. The key tool is that of martingales, which allow explicit computations of expectations
and probabilities with often deep and surprising results.

In particular one obtains well-posedness of SDEs in $\R^d$ of the form
\begin{equation}\label{eq:SDE1}
\d X_t=b(X_t)\d t+\sigma(X_t)\d W_t,
\end{equation}
with $b:\R^d\to\R^d$ and $\sigma:\R^d\to \R^d\otimes\R^d$ smooth coefficients and $(W_t)_{t\geq 0}$ a Brownian motion in $\R^d$. However, in general
$X$ is not better than a \emph{measurable} function of $W$. This fact is rarely mentioned in courses of stochastic calculus, and probabilists seem 
used to it. Nevertheless, a physicist may point out that Brownian motion or its derivative, white noise, are an approximation of a real noise, not
the other way round;  an analyst may found this lack of continuity disturbing. Therefore a theory which is too sensitive on the structure of the noise is not so satisfactory after all. A \emph{robust} theory would be
more convincing from this point of view. 
In the late '70s, the works of Doss \cite{doss77} and Sussmann \cite{sussmann78} gave sufficient conditions on the coefficient 
$\sigma$ for continuity of the maps $W\mapsto X$ in the sup-norm topology on $C([0,T];\R^d)$. These conditions were however very restrictive for $d>1$.

Following an early intuition by Föllmer \cite{follmer81}, Lyons constructed a deterministic (\emph{pathwise}) approach to stochastic integration. The main result is the 
construction of a topology that makes the map $W\mapsto X$ continuous. However, there is a very important twist: the topology is not just on $W$ or $X$, but on a richer object which contains more information. If for example $W:[0,T]\to\R^d$ is a deterministic smooth path, then one needs to consider a finite number
of \emph{iterated integrals} of $W$, which take the form
\[
{\bf W}^{n}_{s,t}=\int_{s<u_1<\cdots<u_n<t} \d{W}_{u_1}\otimes\cdots\otimes\d{W}_{u_n}, \qquad n\in\N, \ 0\leq s\leq t\leq T,
\]
where $\d W_u=\dot{W}_{u}\d u$. For a fixed $\gamma\in\,]0,1[$, one takes $N\in\N$ such that $N\gamma\leq 1<(N+1)\gamma$ and for every
smooth $W:[0,T]\to\R^d$ 
\[
{\bf W}^{(N)}_{s,t}:=1+\sum_{n=1}^N {\bf W}^{n}_{s,t},\qquad 0\leq s\leq t\leq T,
\]
which belongs to the truncated tensor algebra $T^{(N)}=\oplus_{n=0}^N (\R^d)^{\otimes n}$. We note that ${\bf W}^{1}_{s,t}=W_t-W_s$, so that
${\bf W}^{(N)}_{s,t}$ \emph{contains} the increments of the original process, plus additional information. We can now define a distance between
two such objects ${\bf W}^{(N)}_{s,t}$ and ${\bf V}^{(N)}_{s,t}$, for smooth $W,V:[0,T]\to\R^d$
\[
\d_\gamma\left({\bf W}^{(N)},{\bf V}^{(N)}\right):=\sup_{n=1,\ldots,N} \sup_{s\ne t} \frac{\left|{\bf W}^{n}_{s,t}-{\bf V}^{n}_{s,t}\right|}{|t-s|^{n\gamma}}.
\]
Then Lyons' result was that the map ${\bf W}^{(N)}\mapsto {\bf X}^{(N)}$, where $W,X:[0,T]\to\R^d$ are smooth processes which satisfy \eqref{eq:SDE1},
is \emph{continuous} with respect to the metric ${\rm d}_\gamma$.

Lyons' paper \cite{lyons98} was astounding for its novelty: it introduced in stochastic analysis a number of concepts which were unknown to many 
probabilists, in particular the algebraic language based on the work of Chen \cite{chen57} on iterated integrals. Moreover it presented a radically different approach
to the pillar of modern probability theory, the Itô stochastic calculus. For these reasons, it seems that Lyons' ideas took some time before being 
widely accepted by the community and became really famous only fifteen years later, when Hairer proved their power in the context of SPDEs.
See the book of Friz-Hairer \cite{fh14} for a pedagogical introduction.

\subsection{Singular SPDEs and regularity structures}
As we have seen above, several interesting physical models were described in the '80s with SPDEs such as the \emph{dynamical $\phi^4_d$ model},
recall the stochastic quantization \eqref{eq:pw},
\begin{equation}\label{phi4d}
\frac{\partial \phi}{\partial t}= \Delta \phi -\phi^3 + \xi, \qquad x\in\R^d,
\end{equation}
for $d=2,3$ and the KPZ equation \eqref{eq:kpz2}. In both equations there are ill-defined non-linear functionals of some distribution. Equations of
this kind are now commonly known as \emph{singular SPDEs}.

In 2003 Da Prato-Debussche \cite{dpd03} solved the stochastic quantization in $d=2$ with the following idea: they wrote $\phi=z+v$, where
$z$ is the solution to the linear stochastic heat equation with additive white noise
\[
\frac{\partial z}{\partial t}= \Delta z + \xi, \qquad x\in\R^2,
\]
and they wrote an equation for $v=\phi-z$
\[
\frac{\partial v}{\partial t}= \Delta v -z^3-3z^2v-3zv^2-v^3, 
\]
which is now random only through the explicit Gaussian process $z$. We note that $z$ is still a distribution, so that the terms
$z^2$ and $z^3$ are still ill-defined; however it turns out that it is possible to give a meaning to these terms as distributions with the classical 
\emph{Wick renormalization}. Then, the products $z^2v$ and $zv^2$ are defined using \emph{Besov spaces}. This allows to use a fixed point
argument for $v$ and obtain existence and uniqueness for the original (renormalized) equation. 
However this technique does not work for $d=3$, since in this case the products $z^2v$ and $zv^2$ are still ill-defined. 

Since Lyons' foundational paper of 1998, rough 
paths have been based on \emph{generalised Taylor expansions}, with standard monomials replaced by iterated integrals of the driving noise. 
In 2004 Gubinelli built on this idea a new approach to rough integration based on the notion of \emph{controlled paths} \cite{gubi04} and started to work
on the project of a rough approach to SPDEs, see for example the 2010 paper \cite{gt10} with Tindel.

In 2011 Hairer \cite{hairer11} considered the equation
\[
\frac{\partial u}{\partial t} = \frac12\frac{\partial^2 u}{\partial x^2} +g(u)\,\frac{\partial u}{\partial x}+ \xi , \qquad t\geq 0, \ x\in\R
\]
with $u$ and $\xi$ taking values in $\R^d$ with $d>1$, and $g$ taking values in $\R^{d\times d}$.
Although this is less frightening than KPZ, the product $g(u)\frac{\partial u}{\partial x}$ is ill-defined for the usual reason: the partial derivative of $u$
is a distribution, the function $g(u)$ is not smooth, and therefore the product cannot be defined by an integration by parts or other classical tools (the fact
that $u$ is vector valued prevents in general this product from being written as $\frac{\partial}{\partial x}G(u)$).
The idea was to treat the solution $u(t,x)$ as a rough path \emph{in space}. 

In 2013 Hairer managed to apply the same techniques to KPZ \cite{hairer13},
thus giving a well-posedness theory for this equation first introduced in 1986. The importance of this result was amplified by the explosion of
activity around the KPZ universality class following the 2011 papers by Bal\'{a}zs-Quastel-Sepp\"{a}l\"{a}inen \cite{bqs11} and
Amir-Corwin-Quastel \cite{acq11}, which proved that the Cole-Hopf
solution proposed by Bertini-Cancrini has indeed the scaling computed in the original KPZ paper \cite{kpz86} with non-rigorous renormalization group techniques.

In order to solve the stochastic quantization in $d=3$, and many other equations, Hairer \cite{hairer14} expanded the theory of rough paths to cover functions of space-time. Da Prato-Debussche \cite{dpd03} had solved the case $d=2$ with the \emph{global} expansion $\phi=z+v$ of the solution, in terms of an
explicit term $z$ and a \emph{remainder} $v$. Hairer's idea was to use
rather \emph{local} expansions at each point $(t,x)$ in space-time, with a far-reaching generalization of the classical notion of Taylor expansion.
The theory has been developed and expanded in three subsequent papers: Bruned-Hairer-Zambotti \cite{bhz}, Chandra-Hairer \cite{ch16}, 
Bruned-Chandra-Chevyrev-Hairer \cite{BCCH}.

In the meantime, Gubinelli-Imkeller-Perkowski \cite{gip} constructed a different approach to singular SPDEs based on \emph{paracontrolled distributions},
combining the \emph{paradifferential calculus} coming from harmonic analysis and the ideas of rough paths. This approach is effective in many situations like KPZ and the stochastic quantization, see also the paper \cite{mw17} by Mourrat-Weber on the convergence of the two-dimensional dynamic {I}sing-{K}ac model to the dynamical $\phi^4_2$, see \eqref{phi4d},
but not in all cases which are covered by regularity structures. In my personal opinion it is Hairer's theory which
transposes in the most faithful way Gubinelli's ideas on rough paths from SDEs to SPDEs. 

Another interesting approach to the KPZ equation is that of energy solutions by Gon\c{c}alves-Jara \cite{gj14} and Gubinelli-Jara \cite{gj13}, which is 
particularly effective in order to prove convergence under rescaling of a large class of particle systems to a martingale problem formulation of KPZ. 
Uniqueness for such a martingale problem was proved in \cite{gp18} by Gubinelli-Perkowski. Other construction of the $\phi^4_3$ dynamical model are due to Kupiainen \cite{kupiainen}, using renormalization group methods, and to Albeverio-Kusuoka \cite{alku}, using finite-dimensional approximations.

\section{Conclusions}

In this brief and personal history of SPDEs I have left aside many topics that would deserve more attention, for example
\begin{itemize}
\item \emph{regularization by noise}, see Flandoli-Gubinelli-Priola \cite{fgp10}
\item the stochastic FKPP equation, see Mueller-Mytnik-Quastel \cite{mmq11}
\item stochastic dispersive equations, stochastic conservation laws and viscosity solutions for 
fully non-linear SPDEs
\item numerical analysis of SPDEs.
\end{itemize}

I hope that I have at least managed to express my enthousiasm for this topic. The last seven years have been particularly exciting: Gubinelli and Hairer have clearly 
influenced each other in a number of occasions, and their work has spurred an exceptional activity in this area. Rough paths and regularity structures tend to 
make relatively little
use of classical probability theory, and my project of combining stochastic calculus and SPDEs went exactly in the opposite direction. However 
in the years before 2013 I felt somewhat discouraged by the lack of progress of this project, and Hairer's paper on KPZ came as a revelation to me. 
What came afterwards was one of those rare situations when reality surpasses our own dreams.

The message that I wished to convey is that the ground for the success of today was prepared by a considerable amount of work by a whole community, in 
particular on equations driven by space-time white noise. I am convinced that this activity has produced many ideas which could and should be of interest for 
other communities and there are already encouraging signs in this direction.

\bibliographystyle{arxiv}
\bibliography{DCDS}

\providecommand{\href}[2]{#2}\begingroup\raggedright\begin{thebibliography}{10}

\bibitem{ah77}
S.~Albeverio and R.~H{\o}egh-Krohn, ``Dirichlet forms and diffusion processes
  on rigged {H}ilbert spaces,''
  \href{http://dx.doi.org/10.1007/BF00535706}{{\em Z.
  Wahrscheinlichkeitstheorie und Verw. Gebiete} {\bfseries 40} no.~1, (1977)
  1--57}. \url{https://doi.org/10.1007/BF00535706}.

\bibitem{alku}
S.~{Albeverio} and S.~{Kusuoka}, ``{The invariant measure and the flow
  associated to the $\Phi^4_3$-quantum field model},'' {\em ArXiv e-prints}
  (2017) , \href{http://arxiv.org/abs/1711.07108}{{\ttfamily
  arXiv:1711.07108}}.

\bibitem{acq11}
G.~Amir, I.~Corwin, and J.~Quastel, ``Probability distribution of the free
  energy of the continuum directed random polymer in {$1+1$} dimensions,''
  \href{http://dx.doi.org/10.1002/cpa.20347}{{\em Comm. Pure Appl. Math.}
  {\bfseries 64} no.~4, (2011) 466--537}.
  \url{https://doi.org/10.1002/cpa.20347}.

\bibitem{bqs11}
M.~Bal\'{a}zs, J.~Quastel, and T.~Sepp\"{a}l\"{a}inen, ``Fluctuation exponent
  of the {KPZ}/stochastic {B}urgers equation,''
  \href{http://dx.doi.org/10.1090/S0894-0347-2011-00692-9}{{\em J. Amer. Math.
  Soc.} {\bfseries 24} no.~3, (2011) 683--708}.
  \url{https://doi.org/10.1090/S0894-0347-2011-00692-9}.

\bibitem{bt72}
A.~Bensoussan and R.~Temam, ``\'{E}quations aux d\'{e}riv\'{e}es partielles
  stochastiques non lin\'{e}aires. {I},''
  \href{http://dx.doi.org/10.1007/BF02761449}{{\em Israel J. Math.} {\bfseries
  11} (1972) 95--129}. \url{https://doi.org/10.1007/BF02761449}.

\bibitem{bt73}
A.~Bensoussan and R.~Temam, ``\'{E}quations stochastiques du type
  {N}avier-{S}tokes,''
  \href{http://dx.doi.org/10.1016/0022-1236(73)90045-1}{{\em J. Functional
  Analysis} {\bfseries 13} (1973) 195--222}.
  \url{https://doi.org/10.1016/0022-1236(73)90045-1}.

\bibitem{bc95}
L.~Bertini and N.~Cancrini, ``The stochastic heat equation: {F}eynman-{K}ac
  formula and intermittence,'' \href{http://dx.doi.org/10.1007/BF02180136}{{\em
  J. Statist. Phys.} {\bfseries 78} no.~5-6, (1995) 1377--1401}.
  \url{https://doi.org/10.1007/BF02180136}.

\bibitem{bg97}
L.~Bertini and G.~Giacomin, ``Stochastic {B}urgers and {KPZ} equations from
  particle systems,'' \href{http://dx.doi.org/10.1007/s002200050044}{{\em Comm.
  Math. Phys.} {\bfseries 183} no.~3, (1997) 571--607}.
  \url{https://doi.org/10.1007/s002200050044}.

\bibitem{bkl01}
J.~Bricmont, A.~Kupiainen, and R.~Lefevere, ``Ergodicity of the 2{D}
  {N}avier-{S}tokes equations with random forcing,''
  \href{http://dx.doi.org/10.1007/s002200100510}{{\em Comm. Math. Phys.}
  {\bfseries 224} no.~1, (2001) 65--81}.
  \url{https://doi.org/10.1007/s002200100510}. Dedicated to Joel L. Lebowitz.

\bibitem{bkl02}
J.~Bricmont, A.~Kupiainen, and R.~Lefevere, ``Exponential mixing of the 2{D}
  stochastic {N}avier-{S}tokes dynamics,''
  \href{http://dx.doi.org/10.1007/s00220-002-0708-1}{{\em Comm. Math. Phys.}
  {\bfseries 230} no.~1, (2002) 87--132}.
  \url{https://doi.org/10.1007/s00220-002-0708-1}.

\bibitem{BCCH}
Y.~{Bruned}, A.~{Chandra}, I.~{Chevyrev}, and M.~{Hairer}, ``{Renormalising
  SPDEs in regularity structures},'' {\em to appear in J. Eur. Math. Soc.
  (JEMS)} (Nov., 2017) , \href{http://arxiv.org/abs/1711.10239}{{\ttfamily
  arXiv:1711.10239 [math.AP]}}.

\bibitem{bhz}
Y.~Bruned, M.~Hairer, and L.~Zambotti, ``Algebraic renormalisation of
  regularity structures,''
  \href{http://dx.doi.org/10.1007/s00222-018-0841-x}{{\em Invent. Math.}
  {\bfseries 215} no.~3, (2019) 1039--1156}.
  \url{https://doi.org/10.1007/s00222-018-0841-x}.

\bibitem{cabana70}
E.~Caba\~{n}a, ``The vibrating string forced by white noise,''
  \href{http://dx.doi.org/10.1007/BF00531880}{{\em Z.
  Wahrscheinlichkeitstheorie und Verw. Gebiete} {\bfseries 15} (1970)
  111--130}. \url{https://doi.org/10.1007/BF00531880}.

\bibitem{cepa}
E.~C\'{e}pa, ``Probl\`eme de {S}korohod multivoque,''
  \href{http://dx.doi.org/10.1214/aop/1022855642}{{\em Ann. Probab.} {\bfseries
  26} no.~2, (1998) 500--532}. \url{https://doi.org/10.1214/aop/1022855642}.

\bibitem{cerrai01}
S.~Cerrai, \href{http://dx.doi.org/10.1007/b80743}{{\em Second order {PDE}'s in
  finite and infinite dimension}}, vol.~1762 of {\em Lecture Notes in
  Mathematics}.
\newblock Springer-Verlag, Berlin, 2001.
\newblock \url{https://doi.org/10.1007/b80743}.
\newblock A probabilistic approach.

\bibitem{ch16}
A.~Chandra and M.~Hairer, ``An analytic {BPHZ} theorem for {R}egularity
  {S}tructures,'' {\em ArXiv e-prints} (Oct., 2016) ,
  \href{http://arxiv.org/abs/1612.08138}{{\ttfamily arXiv:1612.08138
  [math.AP]}}.

\bibitem{chen57}
K.-T. Chen, ``Integration of paths, geometric invariants and a generalized
  {B}aker-{H}ausdorff formula,'' {\em Ann. of Math. (2)} {\bfseries 65} (1957)
  163--178.

\bibitem{chen64}
Y.~M. Chen, ``On scattering of waves by objects imbedded in random media:
  {S}tochastic linear partial differential equations and scattering of waves by
  conducting sphere imbedded in random media,''
  \href{http://dx.doi.org/10.1063/1.1931186}{{\em J. Mathematical Phys.}
  {\bfseries 5} (1964) 1541--1546}. \url{https://doi.org/10.1063/1.1931186}.

\bibitem{chen15}
Y.-T. Chen, ``Pathwise nonuniqueness for the {SPDE}s of some super-{B}rownian
  motions with immigration,'' \href{http://dx.doi.org/10.1214/14-AOP962}{{\em
  Ann. Probab.} {\bfseries 43} no.~6, (2015) 3359--3467}.
  \url{https://doi.org/10.1214/14-AOP962}.

\bibitem{corwin16}
I.~Corwin, ``Kardar-{P}arisi-{Z}hang universality,''
  \href{http://dx.doi.org/10.1090/noti1334}{{\em Notices Amer. Math. Soc.}
  {\bfseries 63} no.~3, (2016) 230--239}.
  \url{https://doi.org/10.1090/noti1334}.

\bibitem{dpd03}
G.~Da~Prato and A.~Debussche, ``Strong solutions to the stochastic quantization
  equations,'' \href{http://dx.doi.org/10.1214/aop/1068646370}{{\em Ann.
  Probab.} {\bfseries 31} no.~4, (2003) 1900--1916}.
  \url{https://doi.org/10.1214/aop/1068646370}.

\bibitem{dpit76}
G.~Da~Prato, M.~Iannelli, and L.~Tubaro, ``Stochastic differential equations in
  {B}anach spaces, variational formulation,'' {\em Atti Accad. Naz. Lincei
  Rend. Cl. Sci. Fis. Mat. Nat. (8)} {\bfseries 61} no.~3-4, (1976) 168--176
  (1977).

\bibitem{dpz1}
G.~Da~Prato and J.~Zabczyk,
  \href{http://dx.doi.org/10.1017/CBO9780511666223}{{\em Stochastic equations
  in infinite dimensions}}, vol.~44 of {\em Encyclopedia of Mathematics and its
  Applications}.
\newblock Cambridge University Press, Cambridge, 1992.
\newblock \url{https://doi.org/10.1017/CBO9780511666223}.

\bibitem{dpz2}
G.~Da~Prato and J.~Zabczyk,
  \href{http://dx.doi.org/10.1017/CBO9780511662829}{{\em Ergodicity for
  infinite-dimensional systems}}, vol.~229 of {\em London Mathematical Society
  Lecture Note Series}.
\newblock Cambridge University Press, Cambridge, 1996.
\newblock \url{https://doi.org/10.1017/CBO9780511662829}.

\bibitem{dpz3}
G.~Da~Prato and J.~Zabczyk,
  \href{http://dx.doi.org/10.1017/CBO9780511543210}{{\em Second order partial
  differential equations in {H}ilbert spaces}}, vol.~293 of {\em London
  Mathematical Society Lecture Note Series}.
\newblock Cambridge University Press, Cambridge, 2002.
\newblock \url{https://doi.org/10.1017/CBO9780511543210}.

\bibitem{dmz06}
R.~C. Dalang, C.~Mueller, and L.~Zambotti, ``Hitting properties of parabolic
  s.p.d.e.'s with reflection,''
  \href{http://dx.doi.org/10.1214/009117905000000792}{{\em Ann. Probab.}
  {\bfseries 34} no.~4, (2006) 1423--1450}.
  \url{https://doi.org/10.1214/009117905000000792}.

\bibitem{dalecki66}
J.~L. Dalecki\u{\i}, ``Differential equations with functional derivatives and
  stochastic equations for generalized random processes,'' {\em Dokl. Akad.
  Nauk SSSR} {\bfseries 166} (1966) 1035--1038.

\bibitem{dawson72}
D.~A. Dawson, ``Stochastic evolution equations,''
  \href{http://dx.doi.org/10.1016/0025-5564(72)90039-9}{{\em Math. Biosci.}
  {\bfseries 15} (1972) 287--316}.
  \url{https://doi.org/10.1016/0025-5564(72)90039-9}.

\bibitem{dawson93}
D.~A. Dawson, \href{http://dx.doi.org/10.1007/BFb0084190}{``Measure-valued
  {M}arkov processes,''} in {\em \'{E}cole d'\'{E}t\'{e} de {P}robabilit\'{e}s
  de {S}aint-{F}lour {XXI}---1991}, vol.~1541 of {\em Lecture Notes in Math.},
  pp.~1--260.
\newblock Springer, Berlin, 1993.
\newblock \url{https://doi.org/10.1007/BFb0084190}.

\bibitem{dh79}
D.~A. Dawson and K.~J. Hochberg, ``The carrying dimension of a stochastic
  measure diffusion,'' {\em Ann. Probab.} {\bfseries 7} no.~4, (1979) 693--703.
  \url{http://links.jstor.org/sici?sici=0091-1798(197908)7:4<693:TCDOAS>2.0.CO;2-E&origin=MSN}.

\bibitem{doss77}
H.~Doss, ``Liens entre \'{e}quations diff\'{e}rentielles stochastiques et
  ordinaires,'' {\em Ann. Inst. H. Poincar\'{e} Sect. B (N.S.)} {\bfseries 13}
  no.~2, (1977) 99--125.

\bibitem{ems01}
W.~E, J.~C. Mattingly, and Y.~Sinai, ``Gibbsian dynamics and ergodicity for the
  stochastically forced {N}avier-{S}tokes equation,''
  \href{http://dx.doi.org/10.1007/s002201224083}{{\em Comm. Math. Phys.}
  {\bfseries 224} no.~1, (2001) 83--106}.
  \url{https://doi.org/10.1007/s002201224083}. Dedicated to Joel L. Lebowitz.

\bibitem{el94}
K.~D. Elworthy and X.-M. Li, ``Formulae for the derivatives of heat
  semigroups,'' \href{http://dx.doi.org/10.1006/jfan.1994.1124}{{\em J. Funct.
  Anal.} {\bfseries 125} no.~1, (1994) 252--286}.
  \url{https://doi.org/10.1006/jfan.1994.1124}.

\bibitem{el15}
A.~M. Etheridge and C.~Labb\'{e}, ``Scaling limits of weakly asymmetric
  interfaces,'' \href{http://dx.doi.org/10.1007/s00220-014-2243-2}{{\em Comm.
  Math. Phys.} {\bfseries 336} no.~1, (2015) 287--336}.
  \url{https://doi.org/10.1007/s00220-014-2243-2}.

\bibitem{fgp10}
F.~Flandoli, M.~Gubinelli, and E.~Priola, ``Well-posedness of the transport
  equation by stochastic perturbation,''
  \href{http://dx.doi.org/10.1007/s00222-009-0224-4}{{\em Invent. Math.}
  {\bfseries 180} no.~1, (2010) 1--53}.
  \url{https://doi.org/10.1007/s00222-009-0224-4}.

\bibitem{fm95}
F.~Flandoli and B.~Maslowski, ``Ergodicity of the {$2$}-{D} {N}avier-{S}tokes
  equation under random perturbations,'' {\em Comm. Math. Phys.} {\bfseries
  172} no.~1, (1995) 119--141.
  \url{http://projecteuclid.org/euclid.cmp/1104273961}.

\bibitem{follmer81}
H.~F\"{o}llmer, ``Calcul d'{I}t\^{o} sans probabilit\'{e}s,'' in {\em Seminar
  on {P}robability, {XV} ({U}niv. {S}trasbourg, {S}trasbourg, 1979/1980)
  ({F}rench)}, vol.~850 of {\em Lecture Notes in Math.}, pp.~143--150.
\newblock Springer, Berlin, 1981.

\bibitem{fh14}
P.~K. Friz and M.~Hairer,
  \href{http://dx.doi.org/10.1007/978-3-319-08332-2}{{\em A course on rough
  paths}}.
\newblock Universitext. Springer, Cham, 2014.
\newblock \url{https://doi.org/10.1007/978-3-319-08332-2}.
\newblock With an introduction to regularity structures.

\bibitem{fukushima80}
M.~Fukushima, {\em Dirichlet forms and {M}arkov processes}, vol.~23 of {\em
  North-Holland Mathematical Library}.
\newblock North-Holland Publishing Co., Amsterdam-New York; Kodansha, Ltd.,
  Tokyo, 1980.

\bibitem{fot11}
M.~Fukushima, Y.~Oshima, and M.~Takeda, {\em Dirichlet forms and symmetric
  {M}arkov processes}, vol.~19 of {\em De Gruyter Studies in Mathematics}.
\newblock Walter de Gruyter \& Co., Berlin, extended~ed., 2011.

\bibitem{fo01}
T.~Funaki and S.~Olla, ``Fluctuations for {$\nabla\phi$} interface model on a
  wall,'' \href{http://dx.doi.org/10.1016/S0304-4149(00)00104-6}{{\em
  Stochastic Process. Appl.} {\bfseries 94} no.~1, (2001) 1--27}.
  \url{https://doi.org/10.1016/S0304-4149(00)00104-6}.

\bibitem{gibson67}
W.~E. Gibson, ``An exact solution for a class of stochastic partial
  differential equations,'' \href{http://dx.doi.org/10.1137/0115118}{{\em SIAM
  J. Appl. Math.} {\bfseries 15} (1967) 1357--1362}.
  \url{https://doi.org/10.1137/0115118}.

\bibitem{gj87}
J.~Glimm and A.~Jaffe, \href{http://dx.doi.org/10.1007/978-1-4612-4728-9}{{\em
  Quantum physics}}.
\newblock Springer-Verlag, New York, second~ed., 1987.
\newblock \url{https://doi.org/10.1007/978-1-4612-4728-9}.
\newblock A functional integral point of view.

\bibitem{gj14}
P.~Gon\c{c}alves and M.~Jara, ``Nonlinear fluctuations of weakly asymmetric
  interacting particle systems,''
  \href{http://dx.doi.org/10.1007/s00205-013-0693-x}{{\em Arch. Ration. Mech.
  Anal.} {\bfseries 212} no.~2, (2014) 597--644}.
  \url{https://doi.org/10.1007/s00205-013-0693-x}.

\bibitem{gross67}
L.~Gross, ``Potential theory on {H}ilbert space,''
  \href{http://dx.doi.org/10.1016/0022-1236(67)90030-4}{{\em J. Functional
  Analysis} {\bfseries 1} (1967) 123--181}.
  \url{https://doi.org/10.1016/0022-1236(67)90030-4}.

\bibitem{gubi04}
M.~Gubinelli, ``Controlling rough paths,''
  \href{http://dx.doi.org/10.1016/j.jfa.2004.01.002}{{\em Journal of Functional
  Analysis} {\bfseries 216} no.~1, (2004) 86 -- 140}.
  \url{http://www.sciencedirect.com/science/article/pii/S0022123604000497}.

\bibitem{gip}
M.~Gubinelli, P.~Imkeller, and N.~Perkowski, ``Paracontrolled distributions and
  singular {PDE}s,'' \href{http://dx.doi.org/10.1017/fmp.2015.2}{{\em Forum
  Math. Pi} {\bfseries 3} (2015) e6, 75},
  \href{http://arxiv.org/abs/1210.2684}{{\ttfamily 1210.2684}}.
  \url{http://0-dx.doi.org.pugwash.lib.warwick.ac.uk/10.1017/fmp.2015.2}.

\bibitem{gj13}
M.~Gubinelli and M.~Jara, ``Regularization by noise and stochastic {B}urgers
  equations,'' \href{http://dx.doi.org/10.1007/s40072-013-0011-5}{{\em Stoch.
  Partial Differ. Equ. Anal. Comput.} {\bfseries 1} no.~2, (2013) 325--350}.
  \url{https://doi.org/10.1007/s40072-013-0011-5}.

\bibitem{gp18}
M.~Gubinelli and N.~Perkowski, ``Energy solutions of {KPZ} are unique,''
  \href{http://dx.doi.org/10.1090/jams/889}{{\em J. Amer. Math. Soc.}
  {\bfseries 31} no.~2, (2018) 427--471}.
  \url{https://doi.org/10.1090/jams/889}.

\bibitem{gt10}
M.~Gubinelli and S.~Tindel, ``Rough evolution equations,''
  \href{http://dx.doi.org/10.1214/08-AOP437}{{\em Ann. Probab.} {\bfseries 38}
  no.~1, (2010) 1--75}. \url{https://doi.org/10.1214/08-AOP437}.

\bibitem{hairer11}
M.~Hairer, ``Rough stochastic {PDE}s,''
  \href{http://dx.doi.org/10.1002/cpa.20383}{{\em Comm. Pure Appl. Math.}
  {\bfseries 64} no.~11, (2011) 1547--1585}.
  \url{https://doi.org/10.1002/cpa.20383}.

\bibitem{hairer13}
M.~Hairer, ``Solving the {KPZ} equation,''
  \href{http://dx.doi.org/10.4007/annals.2013.178.2.4}{{\em Ann. of Math. (2)}
  {\bfseries 178} no.~2, (2013) 559--664}.
  \url{https://doi.org/10.4007/annals.2013.178.2.4}.

\bibitem{hairer14}
M.~Hairer, ``A theory of regularity structures,''
  \href{http://dx.doi.org/10.1007/s00222-014-0505-4}{{\em Invent. Math.}
  {\bfseries 198} no.~2, (2014) 269--504},
  \href{http://arxiv.org/abs/1303.5113}{{\ttfamily 1303.5113}}.

\bibitem{hm06}
M.~Hairer and J.~C. Mattingly, ``Ergodicity of the 2{D} {N}avier-{S}tokes
  equations with degenerate stochastic forcing,''
  \href{http://dx.doi.org/10.4007/annals.2006.164.993}{{\em Ann. of Math. (2)}
  {\bfseries 164} no.~3, (2006) 993--1032}.
  \url{https://doi.org/10.4007/annals.2006.164.993}.

\bibitem{jlm85}
G.~Jona-Lasinio and P.~K. Mitter, ``On the stochastic quantization of field
  theory,'' {\em Comm. Math. Phys.} {\bfseries 101} no.~3, (1985) 409--436.
  \url{http://projecteuclid.org/euclid.cmp/1104114183}.

\bibitem{kpz86}
M.~Kardar, G.~Parisi, and Y.-C. Zhang, ``Dynamic Scaling of Growning
  Interfaces,'' \href{http://dx.doi.org/10.1103/PhysRevLett.56.889}{{\em Phys.
  Rev. Lett.} {\bfseries 56} no.~9, (1986) 899, 4}.
  \url{https://doi.org/10.1103/PhysRevLett.56.889}.

\bibitem{ks88}
N.~Konno and T.~Shiga, ``Stochastic partial differential equations for some
  measure-valued diffusions,'' \href{http://dx.doi.org/10.1007/BF00320919}{{\em
  Probab. Theory Related Fields} {\bfseries 79} no.~2, (1988) 201--225}.
  \url{https://doi.org/10.1007/BF00320919}.

\bibitem{krylov94}
N.~V. Krylov, ``A {$W^n_2$}-theory of the {D}irichlet problem for {SPDE}s in
  general smooth domains,'' \href{http://dx.doi.org/10.1007/BF01192260}{{\em
  Probab. Theory Related Fields} {\bfseries 98} no.~3, (1994) 389--421}.
  \url{https://doi.org/10.1007/BF01192260}.

\bibitem{kr77}
N.~V. Krylov and B.~L. Rozovski\u{\i}, ``The {C}auchy problem for linear
  stochastic partial differential equations,'' {\em Izv. Akad. Nauk SSSR Ser.
  Mat.} {\bfseries 41} no.~6, (1977) 1329--1347, 1448.

\bibitem{ks00}
S.~Kuksin and A.~Shirikyan, ``Stochastic dissipative {PDE}s and {G}ibbs
  measures,'' \href{http://dx.doi.org/10.1007/s002200000237}{{\em Comm. Math.
  Phys.} {\bfseries 213} no.~2, (2000) 291--330}.
  \url{https://doi.org/10.1007/s002200000237}.

\bibitem{ks01}
S.~Kuksin and A.~Shirikyan, ``Ergodicity for the randomly forced 2{D}
  {N}avier-{S}tokes equations,''
  \href{http://dx.doi.org/10.1023/A:1011989910997}{{\em Math. Phys. Anal.
  Geom.} {\bfseries 4} no.~2, (2001) 147--195}.
  \url{https://doi.org/10.1023/A:1011989910997}.

\bibitem{kupiainen}
A.~Kupiainen, ``Renormalization group and stochastic {PDE}s,''
  \href{http://dx.doi.org/10.1007/s00023-015-0408-y}{{\em Ann. Henri
  Poincar\'{e}} {\bfseries 17} no.~3, (2016) 497--535}.
  \url{https://doi.org/10.1007/s00023-015-0408-y}.

\bibitem{legall99}
J.-F. Le~Gall, \href{http://dx.doi.org/10.1007/978-3-0348-8683-3}{{\em Spatial
  branching processes, random snakes and partial differential equations}}.
\newblock Lectures in Mathematics ETH Z\"{u}rich. Birkh\"{a}user Verlag, Basel,
  1999.
\newblock \url{https://doi.org/10.1007/978-3-0348-8683-3}.

\bibitem{lyon60}
R.~H. Lyon, ``Response of a nonlinear string to random excitation,''
  \href{http://dx.doi.org/10.1121/1.1908341}{{\em J. Acoust. Soc. Amer.}
  {\bfseries 32} (1960) 953--960}. \url{https://doi.org/10.1121/1.1908341}.

\bibitem{lyons98}
T.~J. Lyons, ``Differential equations driven by rough signals,''
  \href{http://dx.doi.org/10.4171/RMI/240}{{\em Rev. Mat. Iberoamericana}
  {\bfseries 14} no.~2, (1998) 215--310}.
  \url{https://doi.org/10.4171/RMI/240}.

\bibitem{mr92}
Z.~M. Ma and M.~R\"{o}ckner,
  \href{http://dx.doi.org/10.1007/978-3-642-77739-4}{{\em Introduction to the
  theory of (nonsymmetric) {D}irichlet forms}}.
\newblock Universitext. Springer-Verlag, Berlin, 1992.
\newblock \url{https://doi.org/10.1007/978-3-642-77739-4}.

\bibitem{matt99}
J.~C. Mattingly, ``Ergodicity of {$2$}{D} {N}avier-{S}tokes equations with
  random forcing and large viscosity,''
  \href{http://dx.doi.org/10.1007/s002200050706}{{\em Comm. Math. Phys.}
  {\bfseries 206} no.~2, (1999) 273--288}.
  \url{https://doi.org/10.1007/s002200050706}.

\bibitem{miermont}
G.~Miermont, ``Aspects of random maps,''.
  \url{http://perso.ens-lyon.fr/gregory.miermont/coursSaint-Flour.pdf}.
  Saint-Flour Lecture notes, 2014.

\bibitem{mw17}
J.-C. Mourrat and H.~Weber, ``Convergence of the two-dimensional dynamic
  {I}sing-{K}ac model to {$\Phi^4_2$},''
  \href{http://dx.doi.org/10.1002/cpa.21655}{{\em Comm. Pure Appl. Math.}
  {\bfseries 70} no.~4, (2017) 717--812}.
  \url{https://doi.org/10.1002/cpa.21655}.

\bibitem{mueller91}
C.~Mueller, ``On the support of solutions to the heat equation with noise,''
  \href{http://dx.doi.org/10.1080/17442509108833738}{{\em Stochastics
  Stochastics Rep.} {\bfseries 37} no.~4, (1991) 225--245}.
  \url{https://doi.org/10.1080/17442509108833738}.

\bibitem{mmp14}
C.~Mueller, L.~Mytnik, and E.~Perkins, ``Nonuniqueness for a parabolic {SPDE}
  with {$\frac{3}{4}-\epsilon$}-{H}\"{o}lder diffusion coefficients,''
  \href{http://dx.doi.org/10.1214/13-AOP870}{{\em Ann. Probab.} {\bfseries 42}
  no.~5, (2014) 2032--2112}. \url{https://doi.org/10.1214/13-AOP870}.

\bibitem{mmq11}
C.~Mueller, L.~Mytnik, and J.~Quastel, ``Effect of noise on front propagation
  in reaction-diffusion equations of {KPP} type,''
  \href{http://dx.doi.org/10.1007/s00222-010-0292-5}{{\em Invent. Math.}
  {\bfseries 184} no.~2, (2011) 405--453}.
  \url{https://doi.org/10.1007/s00222-010-0292-5}.

\bibitem{mytnik96}
L.~Mytnik, ``Superprocesses in random environments,''
  \href{http://dx.doi.org/10.1214/aop/1041903212}{{\em Ann. Probab.} {\bfseries
  24} no.~4, (1996) 1953--1978}. \url{https://doi.org/10.1214/aop/1041903212}.

\bibitem{mp11}
L.~Mytnik and E.~Perkins, ``Pathwise uniqueness for stochastic heat equations
  with {H}\"{o}lder continuous coefficients: the white noise case,''
  \href{http://dx.doi.org/10.1007/s00440-009-0241-7}{{\em Probab. Theory
  Related Fields} {\bfseries 149} no.~1-2, (2011) 1--96}.
  \url{https://doi.org/10.1007/s00440-009-0241-7}.

\bibitem{nualart}
D.~Nualart, {\em The {M}alliavin calculus and related topics}.
\newblock Probability and its Applications (New York). Springer-Verlag, Berlin,
  second~ed., 2006.

\bibitem{nupa}
D.~Nualart and E.~Pardoux, ``White noise driven quasilinear {SPDE}s with
  reflection,'' \href{http://dx.doi.org/10.1007/BF01195389}{{\em Probab. Theory
  Related Fields} {\bfseries 93} no.~1, (1992) 77--89}.
  \url{https://doi.org/10.1007/BF01195389}.

\bibitem{pardoux72}
E.~Pardoux, ``Sur des \'{e}quations aux d\'{e}riv\'{e}es partielles
  stochastiques monotones,'' {\em C. R. Acad. Sci. Paris S\'{e}r. A-B}
  {\bfseries 275} (1972) A101--A103.

\bibitem{pw81}
G.~Parisi and Y.~S. Wu, ``Perturbation theory without gauge fixing,'' {\em Sci.
  Sinica} {\bfseries 24} no.~4, (1981) 483--496.

\bibitem{perkins02}
E.~Perkins, ``Dawson-{W}atanabe superprocesses and measure-valued diffusions,''
  in {\em Lectures on probability theory and statistics ({S}aint-{F}lour,
  1999)}, vol.~1781 of {\em Lecture Notes in Math.}, pp.~125--324.
\newblock Springer, Berlin, 2002.

\bibitem{quastel12}
J.~Quastel, ``Introduction to {KPZ},'' in {\em Current developments in
  mathematics, 2011}, pp.~125--194.
\newblock Int. Press, Somerville, MA, 2012.

\bibitem{reyo}
D.~Revuz and M.~Yor, \href{http://dx.doi.org/10.1007/978-3-662-06400-9}{{\em
  Continuous martingales and {B}rownian motion}}, vol.~293 of {\em Grundlehren
  der Mathematischen Wissenschaften [Fundamental Principles of Mathematical
  Sciences]}.
\newblock Springer-Verlag, Berlin, third~ed., 1999.
\newblock \url{https://doi.org/10.1007/978-3-662-06400-9}.

\bibitem{rozovski74}
B.~L. Rozovski\u{\i}, ``Stochastic differential equations in
  infinite-dimensional spaces, and filtering problems,'' in {\em Proceedings of
  the {S}chool and {S}eminar on the {T}heory of {R}andom {P}rocesses
  ({D}ruskininkai, 1974), {P}art {II} ({R}ussian)}, pp.~147--194.
\newblock 1975.

\bibitem{rozovski75}
B.~L. Rozovski\u{\i}, ``Stochastic partial differential equations,'' {\em Mat.
  Sb. (N.S.)} {\bfseries 96(138)} (1975) 314--341, 344.

\bibitem{shiga81}
T.~Shiga, ``Diffusion processes in population genetics,''
  \href{http://dx.doi.org/10.1215/kjm/1250522109}{{\em J. Math. Kyoto Univ.}
  {\bfseries 21} no.~1, (1981) 133--151}.
  \url{https://doi.org/10.1215/kjm/1250522109}.

\bibitem{shiga87}
T.~Shiga, ``Existence and uniqueness of solutions for a class of nonlinear
  diffusion equations,'' \href{http://dx.doi.org/10.1215/kjm/1250520714}{{\em
  J. Math. Kyoto Univ.} {\bfseries 27} no.~2, (1987) 195--215}.
  \url{https://doi.org/10.1215/kjm/1250520714}.

\bibitem{simon74}
B.~Simon, {\em The {$P(\phi )_{2}$} {E}uclidean (quantum) field theory}.
\newblock Princeton University Press, Princeton, N.J., 1974.
\newblock Princeton Series in Physics.

\bibitem{spiegel52}
M.~R. Spiegel, ``The random vibrations of a string,''
  \href{http://dx.doi.org/10.1090/qam/45976}{{\em Quart. Appl. Math.}
  {\bfseries 10} (1952) 25--33}. \url{https://doi.org/10.1090/qam/45976}.

\bibitem{sv79}
D.~W. Stroock and S.~R.~S. Varadhan, {\em Multidimensional diffusion
  processes}.
\newblock Classics in Mathematics. Springer-Verlag, Berlin, 2006.
\newblock Reprint of the 1997 edition.

\bibitem{sussmann78}
H.~J. Sussmann, ``On the gap between deterministic and stochastic ordinary
  differential equations,''
  \href{http://dx.doi.org/10.1016/0166-218x(83)90112-9}{{\em Ann. Probability}
  {\bfseries 6} no.~1, (1978) 19--41}.
  \url{https://doi.org/10.1016/0166-218x(83)90112-9}.

\bibitem{walsh86}
J.~B. Walsh, \href{http://dx.doi.org/10.1007/BFb0074920}{``An introduction to
  stochastic partial differential equations,''} in {\em \'{E}cole d'\'{e}t\'{e}
  de probabilit\'{e}s de {S}aint-{F}lour, {XIV}---1984}, vol.~1180 of {\em
  Lecture Notes in Math.}, pp.~265--439.
\newblock Springer, Berlin, 1986.
\newblock \url{https://doi.org/10.1007/BFb0074920}.

\bibitem{watanabe68}
S.~Watanabe, ``A limit theorem of branching processes and continuous state
  branching processes,'' \href{http://dx.doi.org/10.1215/kjm/1250524180}{{\em
  J. Math. Kyoto Univ.} {\bfseries 8} (1968) 141--167}.
  \url{https://doi.org/10.1215/kjm/1250524180}.

\bibitem{zakai69}
M.~Zakai, ``On the optimal filtering of diffusion processes,''
  \href{http://dx.doi.org/10.1007/BF00536382}{{\em Z.
  Wahrscheinlichkeitstheorie und Verw. Gebiete} {\bfseries 11} (1969)
  230--243}. \url{https://doi.org/10.1007/BF00536382}.

\bibitem{lz01}
L.~Zambotti, ``A reflected stochastic heat equation as symmetric dynamics with
  respect to the 3-d {B}essel bridge,''
  \href{http://dx.doi.org/10.1006/jfan.2000.3685}{{\em J. Funct. Anal.}
  {\bfseries 180} no.~1, (2001) 195--209}.
  \url{https://doi.org/10.1006/jfan.2000.3685}.

\bibitem{lz02}
L.~Zambotti, ``Integration by parts formulae on convex sets of paths and
  applications to {SPDE}s with reflection,''
  \href{http://dx.doi.org/10.1007/s004400200203}{{\em Probab. Theory Related
  Fields} {\bfseries 123} no.~4, (2002) 579--600}.
  \url{https://doi.org/10.1007/s004400200203}.

\bibitem{lz04}
L.~Zambotti, ``Occupation densities for {SPDE}s with reflection,''
  \href{http://dx.doi.org/10.1214/aop/1078415833}{{\em Ann. Probab.} {\bfseries
  32} no.~1A, (2004) 191--215}. \url{https://doi.org/10.1214/aop/1078415833}.

\bibitem{lz15}
L.~Zambotti, \href{http://dx.doi.org/10.1007/978-3-319-52096-4}{{\em Random
  obstacle problems}}, vol.~2181 of {\em Lecture Notes in Mathematics}.
\newblock Springer, Cham, 2017.
\newblock \url{https://doi.org/10.1007/978-3-319-52096-4}.
\newblock Lecture notes from the 45th Probability Summer School held in
  Saint-Flour, 2015.

\end{thebibliography}\endgroup

\end{document}